\documentclass[a4paper,11pt]{article}

\usepackage{textcomp}
\usepackage{amsfonts}
\usepackage[dvips]{graphicx}
\usepackage{amsmath}
\usepackage{amssymb}
\usepackage{mathrsfs}

\title{\textbf{Inversely Unstable
Solutions of Two-Dimensional Systems on Genus-\textit{p} Surfaces
and the Topology of Knotted Attractors}}

\author{Yi SONG and Stephen P. BANKS,\\
Department of Automatic Control and Systems Engineering,\\
University of Sheffield, Mappin Street,\\
Sheffield, S1 3JD.\\
e-mail: s.banks@sheffield.ac.uk}

\begin{document}

\maketitle

\newtheorem{theorem}{Theorem}[section]
\newtheorem{lemma}{Lemma}[section]
\newtheorem{definition}{Definition}[section]
\newtheorem{example}{Example}[section]

\begin{abstract}

In this paper, we will show that a periodic nonlinear ,time-varying
dissipative system that is defined on a genus-$p$ surface contains
one or more invariant sets which act as attractors. Moreover, we
shall generalize a result in [Martins, 2004] and give conditions
under which these invariant sets are not homeomorphic to a circle
individually, which implies the existence of chaotic behaviour. This
is achieved by studying the appearance
of inversely unstable solutions within each invariant set.
\end{abstract}

\section{Introduction}

The general theory of dynamical systems is, of course, a subject
with a long and distinguished history , (see, for example, [Smale,
1967], [Bowen, 1928] and [Manning, 1974]). In particular, the study
of the dynamical behaviour of non-conservative and chaotic systems
has attracted a lot of attention in the past, (see, for example,
[Levinson, 1944], [Martins, 2004] and [Wiggins, 1988]). Consider a
system

\begin{equation} \label{general system}
\left\{ \begin{array}{l}
\dot{x}=F(x,y,t)\\
\dot{y}=G(x,y,t),
\end{array}
\right.
\end{equation}

\noindent where $F(x,y,t)$ and $G(x,y,t)$ are both periodic in t. It
will be called dissipative or non-conservative if there is a locally
proper invariant set on the corresponding $2$-manifold on which the
system is defined. Most real systems are of this kind. Up to the
present a great deal of interest has been paid to the study of the
topology of this invariant set (e.g. [Levinson, 1944]).

Recently, in [Martins, 2004], it is shown that a system

\begin{equation}
\ddot{x} +h(x)\dot{x} +g(t,x)=0,
\end{equation}

\noindent where $h$ and $g$ are smooth functions, periodic on both
$x$ and $t$, which is essentially a periodic nonlinear
2-dimensional, time-varying oscillator with appropriate damping,
contains an invariant set which is not homeomorphic to a circle if
there exists an inversely unstable solution.

In this paper we are interested in generalizing this result and we
will show that instead of just one invariant set, several attractors
can coexist and even be linked in a higher genus surface on which
the system is defined. We will also study the topology of these
attractors in a similar way to [Martins, 2004] and show the
existence of an inversely unstable solution implies that the
specific invariant set is not homeomorphic to a circle.

Moreover, in [Banks, 2002], a way to express a system situated on a
genus-$p$ surface in terms of a spherical one is presented. This is
achieved by opening each handle, i.e., cutting along a fundamental
circuit which contains no equilibria and adding appropriate fixed
points on the resulting sphere \big(as shown in fig (\ref{reduce
genus})\big). In this paper, we will also study the relation between
dissipative systems on a $p$-hole surface and those sitting on a
sphere.

\begin{figure}[!hbp]
\begin{center}
\includegraphics[width=4in]{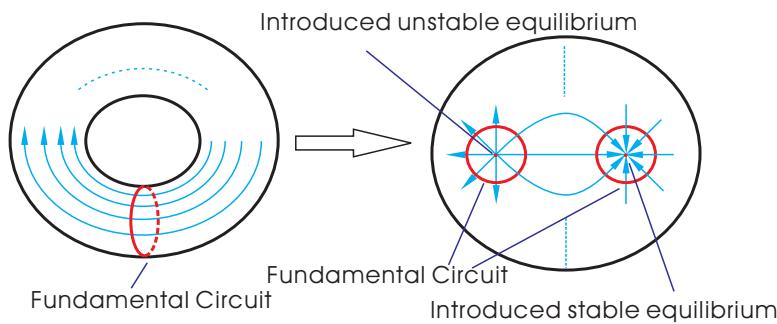}
\caption{Express a Genus-1 System onto a Sphere} \label{reduce
genus}
\end{center}
\end{figure}

In order to motivate the ideas, we reformulate Martins' result in
the following way.

The system given by (\ref{general system}) can be written in the
form

\begin{equation}
\left\{
\begin{array}{l}
\dot{y_1}=y_2-H(y_1)\\
\dot{y_2}=-g(t,y_1)
\end{array}\right.
\end{equation}

\noindent where $H(x)={\int}_0^xh(s)ds$, and $g$ is $T$-periodic in
$t$. The $Poincar\acute{e}$ map is defined as $P(y_0)=y(T;0,y_0)$.
Since the vector field $(y_1,y_2) \to \big( y_2-H(y_1), -g(t,y_1)
\big)$ is periodic with period $R=\big( 1,h(1) \big)$, the solutions
$y$ and $y+kR$ $(k \in \mathbb{Z})$ are equivalent and so the system
may be defined on a cylinder, as in fig (\ref{invariant set on a
cylinder}).

\begin{figure}[!h]
\begin{center}
\includegraphics[width=2in,height=2in]{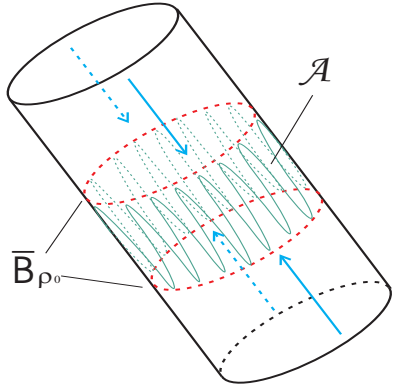}
\caption{The Invariant Set Defined on a Cylinder} \label{invariant
set on a cylinder}
\end{center}
\end{figure}

\noindent Here $\mathcal{A}$ is the invariant set

\[
\mathcal{A} = \bigcap_{n\in
\mathbb{N}}\overline{P}^n(\overline{B}_{\rho_0})
\]

\noindent where $B_{\rho_0}$ is some bounded set \big(which exists
because the system is dissipative, as implied by the arrows in fig
(\ref{invariant set on a cylinder})\big). In [Martins, 2004], he
shows that $\mathcal{A}$ is not homeomorphic to a circle if there is
an inversely unstable periodic orbit somewhere; we can think of the
problem as sitting on a torus with one unstable cycle, as in fig
(\ref{invariant set on torus}). It is in this form that we shall
generalize the result to higher genus surfaces.

\begin{figure}[!h]
\begin{center}
\includegraphics[width=2in,height=2in]{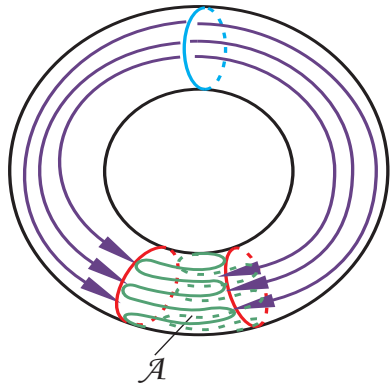}
\caption{Invariant Set in the Torus Case} \label{invariant set on
torus}
\end{center}
\end{figure}

\section{Systems On Genus-p Surfaces}

In [Banks \& Song, 2006] we have shown how to write down analytic
(or meromorphic) systems on genus-$p$ surfaces by the use of
automorphic functions. These systems are not general enough,
however,  to include systems with knots, chaotic annuli, etc.. So we
must consider vector fields which are $C^{\infty}$ but which are
invariant under certain linear, fractional transforms. This will be
the analogue of systems which are periodic in [Martins, 2004] and
have inversely unstable periodic motions -- the latter now becoming
knots on the genus-$p$ surface.

\begin{figure}[!hbp]
\begin{center}
\includegraphics[width=3.5in]{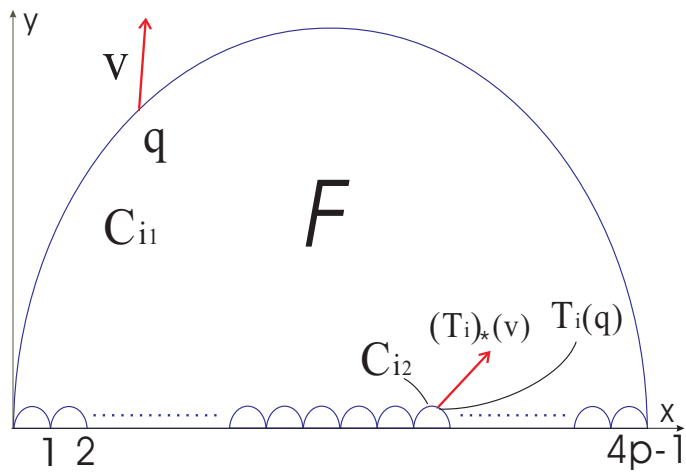}
\caption{Fundamental region of a genus-$p$ surface}
\label{fundamental region}
\end{center}
\end{figure}

In order to generate the most general $C^{\infty}$ systems on
genus-$p$ surfaces, consider a fundamental domain $F$ in the
upper-half plane model of the hyperbolic plane for the surface as in
fig (\ref{fundamental region}).

\noindent Let $\Gamma$ be a $Fuchsian$ group (see [Ford, 1929] and
[Banks \& Song, 2006]) with a subset $\Gamma_1=\{T_i\}$ that
contains maps which pair the sides of $F$. Each map $T_i$ is of the
form:

\begin{equation}
T_i(z)=\frac{az+b}{cz+d}
\end{equation}

\noindent and we shall consider them in real form:

\begin{equation}
T_i(x,y)=\big(\tau_{ix}(x,y), \tau_{iy}(x,y)\big)
\end{equation}

\noindent where $z=x+iy$.

Because of the need to generate $C^{\infty}$ systems defined on a
genus-$p$ surface, we must ensure that, if $T_i$ pairs the sides
$C_{i_1}$ and $C_{i_2}$, as in fig (\ref{fundamental region}), then
the vector field $v$ in the hyperbolic plane at corresponding points
$q$ satisfies

\begin{equation}
v\big(T_i(q)\big)=(T_i)_{\ast}\big(v(q)\big)
\end{equation}

\noindent where $(T_i)_{\ast}$ is the tangent map of $T_i$.

\begin{lemma}
There exists a map from F onto a rectangle R which is one-to-one on
the interior and $C^{\infty}$ $apart$ from at the cusp points.
\end{lemma}

\noindent \textbf{Proof.} We shall construct the map explicitly so
that the required properties will be clear.

\begin{figure}[!hbp]
\begin{center}
\includegraphics[width=4.5in]{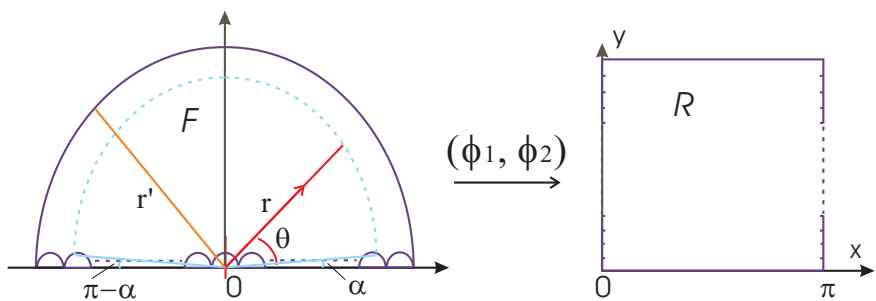}
\caption{Map the fundamental region \textit{F} onto a rectangle
\textit{R}} \label{map onto a rectangle}
\end{center}
\end{figure}

An elementary calculation shows that

\[
\left\{ \begin{array}{l}
x={\phi}_1(r,\theta)\\
y={\phi}_2(r,\theta)
\end{array} \right.
\]
\noindent where
\begin{eqnarray} \label{phi function}
{\phi}_1(r,\theta)&=&\frac{\pi}{\pi-2\alpha}\cdot (\theta-\alpha) \nonumber \\
{\phi}_2(r,\theta)&=&r
\end{eqnarray}

\noindent where $\alpha$ is the value of the starting angle
corresponding to the curve within the fundamental region in the
$(r,\theta)-plane$ \big(as shown in fig (\ref{map onto a
rectangle})\big). \qquad $\Box$

We shall call $R$ the modified fundamental region, and write this
map as

\begin{equation}
\phi: (r,\theta) \to (x,y).
\end{equation}

\noindent Let

\begin{equation}
D_i={\phi}(C_i)
\end{equation}

\noindent be the edges of the curves $C_i$ on the boundary of $F$.
From the above remarks we see that a vector field \textit{w} on $R$
which is associated with one on the original surface must satisfy

\begin{equation} \label{modified system}
w\Big({\phi}\big(T_i(q)\big)\Big)={\phi}_\ast\Big((T_i)_\ast\big(w(q)\big)\Big)
, \quad q\in D_{i_1}
\end{equation}

\noindent where $T_i$ pairs the segments $D_{i_1}$ and $D_{i_2}$.
Let

\[
m_1, m_2, \cdots, m_{4p} \in R,
\]

\noindent denote the points
\[
m_i=\phi(i)
\]

\noindent (i.e., the $\phi$-image of the cusp points on $F$). Then
we have

\begin{lemma}
Any vector field w which is $C^{\infty}$ on the interior of R and
satisfies (\ref{modified system}) where $\phi$ is given by (\ref{phi
function}) and such that
\[
w(m_i)=0
\]
\noindent defines a unique vector field on a genus-p $(p>1)$
surface.
\end{lemma}

\noindent \textbf{Proof.} The only part left to prove is the
converse. This follows from the above remarks and the
${Poincar\acute{e}}$ index theorem---any dynamical system on a
surface of genus $p>1$ must have at least one equilibrium point. We
can choose this as the cusp points of $F$. \qquad $\Box$\\

\begin{figure}[!hbp]
\begin{tabular}{cc}
\includegraphics[width=2in,height=3in]{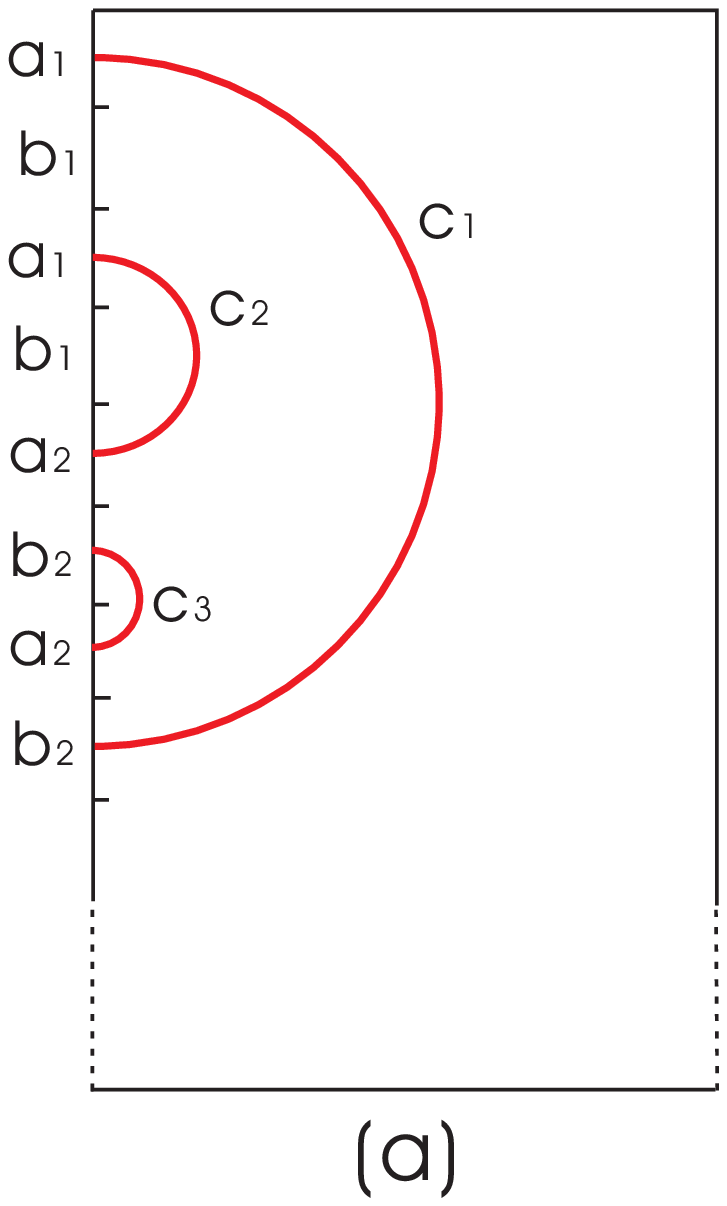}
\includegraphics[width=3in,height=3in]{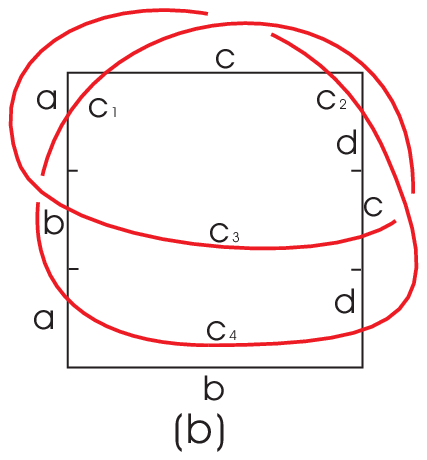}
\end{tabular}
\caption{Closed curves on a surface} \label{knot in a rectangle}
\end{figure}

We next consider the existence of periodic knotted trajectories on
the surface. By the above results we can restrict attention to a
rectangle $R$ as shown in fig (\ref{map onto a rectangle}). Any
closed curve on the surface is given by a set of non-intersecting
curves which `match' in the sense of (\ref{modified system}) at
identified boundaries. For example, the set of curves shown in fig
(\ref{knot in a rectangle}.a) form a closed curve on the
corresponding surface; moreover, fig (\ref{knot in a rectangle}.b)
stands for a \textit{trefoil knot} on a $2$-hole torus if we
identify the sides properly.

\begin{figure}[!hbp]
\begin{center}
  \includegraphics[width=3in]{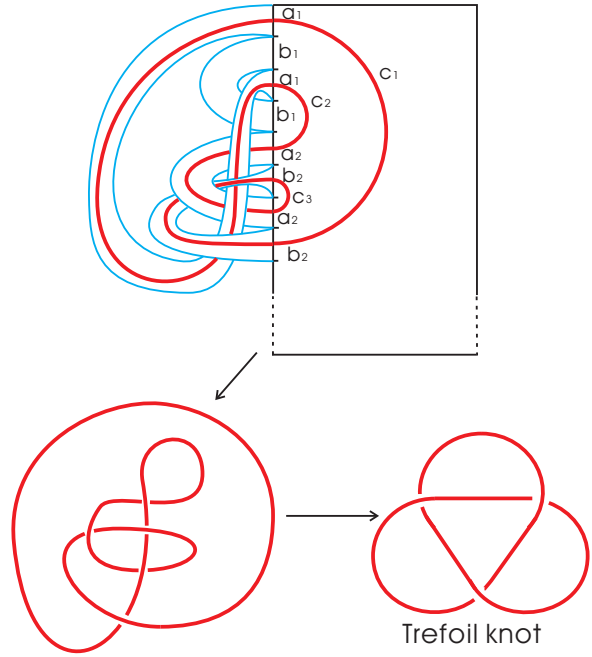}
  \caption{A $Trefoil$ Knot On a Surface}\label{embedding}
  \end{center}
\end{figure}

Of course, the knot type of this closed curve depends on the
embedding of the surface in $\mathbb{R}^3$ (or $\mathrm{S}^3$). For
instance, the surface in fig (\ref{knot in a rectangle}.a) could be
embedded as in fig (\ref{embedding}), which also gives a knot
diagram from which one can calculate a knot invariant (such as the
Kauffman Polynomial).

Let $\psi_i(x,y,t)$ denote the curve $C_i$ within the modified
fundamental region $R$, $f_i(x,y)$ and $g_i(x,y)$ be any
$C^{\infty}$ functions that guarantee the matching of vector fields
at the identified boundaries via (\ref{modified system}). Hence we
have

\begin{lemma}
If there are $C_i$ ($1 \leqslant i \leqslant k$) curves within the
modified fundamental region that stand for a periodic trajectory of
a dynamical system on a genus-p surface in $\mathbb{R}^3$, then this
system can be defined by
\begin{eqnarray}
\dot{x}&=&\sum_{i=1}^{k} \Bigg( \bigg(\frac{\partial\psi_i}{\partial
y}+\psi_if_i \bigg) \cdot \prod_{j \ne i}\psi_j \Bigg) \nonumber\\
\dot{y}&=&\sum_{i=1}^{k} \Bigg(
\bigg(-\frac{\partial\psi_i}{\partial x}+\psi_i g_i \bigg) \cdot
\prod_{j \ne i}\psi_j \Bigg) \label{periodic orbit}
\end{eqnarray}
\end{lemma}

\noindent \textbf{Proof.} Since $\psi_i$ defines the curve $C_i$ in
region $R$, we get
\begin{eqnarray}
\psi_i(x,y,t_i) &=& 0 \label{zero} \\
\qquad \frac{\partial \psi_i}{\partial x} \cdot \dot{x} +
\frac{\partial \psi_i}{\partial y} \cdot \dot{y} &=& 0
\label{condition}
\end{eqnarray}

\noindent For each curve $C_i$, $\psi_i$ switches off all terms in
(\ref{periodic orbit}) except the $i$th one. Substitute
(\ref{periodic orbit}) into (\ref{condition}) and we have
\begin{eqnarray}
{\frac{\partial \psi_i}{\partial x}} \bigg({\frac{\partial
\psi_i}{\partial y}} + \psi_i f_i \bigg){\prod_{j \ne i} {\psi_j}}
+{\frac{\partial \psi_i}{\partial y}} \bigg(-{\frac{\partial
\psi_i}{\partial x}} + \psi_i g_i \bigg){\prod_{j \ne i} {\psi_j}}
\nonumber \\
= {\prod_{j \ne i}} {\psi_j}\bigg({\frac{\partial \psi_i}{\partial
x}}\psi_i f_i+{\frac{\partial \psi_i}{\partial y}}\psi_i
g_i\bigg)=0,
\end{eqnarray}
\noindent so the lemma follows. \qquad $\Box$

\section{The $\mathbf{Poincar\acute{e}}$ Map and Knotted Attractors}

Equation (\ref{periodic orbit}) now can be regarded as a general
form of dynamical systems in the hyperbolic upper-half plane, which
can be situated on a genus-$p$ surface after identification of the
corresponding sides.

Again Consider the ${Poincar\acute{e}}$ map $P: \mathbb{R}^2 \to
\mathbb{R}^2$ given by
\[
P(x_0,y_0)=Y\big(T;0,(x_0,y_0)\big),
\]
\noindent where $Y(t)=Y\big( t;0,(x_0,y_0)\big)$ is the solution of
(\ref{periodic orbit}) starting from point $(x_0,y_0)$. Because of
the `periodicity' from the automorphic form which is defined by the
$Fuchsian$ group $\Gamma$, we have
\begin{equation}
P\big( \Gamma_i(x_0,y_0) \big) = \Gamma_i \big( P(x_0,y_0)
\big)
\end{equation}
\noindent where $\Gamma_i$ is the transformation from one
fundamental region to another one next to it. Moreover, if $(x,y)$
is a solution of (\ref{periodic orbit}), so is $\Gamma_{i}^{n}(x,y)$
$(n \in \mathbb{Z})$. Restricting our attention onto one fundamental
region $F$ \big(as shown in fig (\ref{fundamental region})\big), we
only need to consider the dynamics within it, and obviously the
$Poincar\acute{e}$ map is well defined on $F$.

If the system given by (\ref{periodic orbit}) is dissipative, then
there exists an unstable periodic orbit, which means all the
trajectories are pointing outward along this closed curve. Assume
that it is represented by $\{ \psi_i \} $ $ (1 \le i \le k)$, we are
mainly interested in what the dynamics will look like on the rest of
the surface.

To begin, we need to perform some surgery to the $2$-manifold. By
cutting along this closed orbit, the genus-$p$ surface will
effectively turn into
a $(p-1)$-hole torus with two boundary circles being introduced.\\

\noindent \textit{Remark.} To make this statement much clearer, we
now look at an example of cutting along a \textit{trefoil knot} that
sits on a torus.

\begin{figure}[!h]
\begin{center}
\includegraphics[width=4in]{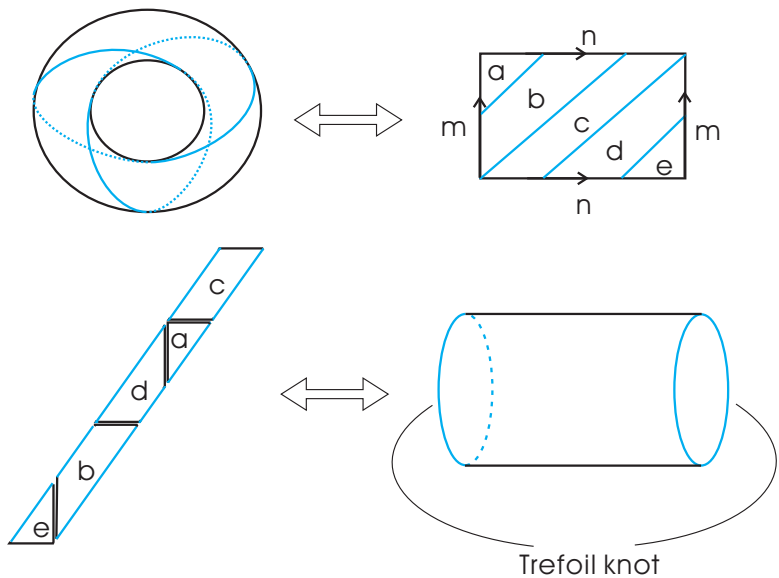}
\caption{Cut a $Trefoil$ Knot On a Torus} \label{cut trefoil}
\end{center}
\end{figure}

As shown in fig (\ref{cut trefoil}), by cutting the torus along this
\textit{trefoil knot} and identifying the corresponding parts on
both sides, `m' and `n', we get a cylinder with the two ends being
the original \textit{trefoil knot}. This surgery can always be
performed on the genus-$p$ surface such that the knot along which is
cut will generate one pair of the sides in the fundamental domain,
and this results in the constructed 2-manifold being a $(g-1)$-hole
torus with two boundary circles \big(as shown in fig (\ref{cut
g-hole torus})\big).

\begin{figure}[!h]
\begin{center}
\includegraphics[width=4in]{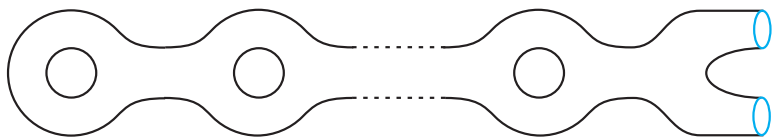}
\caption{Cut a $g$-hole Torus Open Along a Knot} \label{cut g-hole
torus}
\end{center}
\end{figure}

In [Martins, 2004], he studied the torus case, (i.e.,
$\textrm{genus}=1$), and showed that if there exists a trivial
unstable periodic orbit, then an invariant set $\mathcal{A}$, a band
around the tube, which may or may not be homeomorphic to a circle,
must exist (see fig (\ref{invariant set on torus}) for an
illustration). $\mathcal{A}$ is a compact, non-empty, connected set,
and it acts as an attractor towards which all the dynamics converge.
It is given by the iterations of the $Poinca\acute{e}$ map within
the fundamental region to a well-defined bounded set.

\begin{figure}[!h]
\begin{center}
\includegraphics[width=1in]{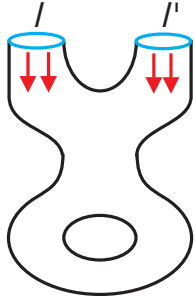}
\caption{Cut a $2$-hole Torus} \label{cut 2-hole torus}
\end{center}
\end{figure}

In the case of genus-$2$ surfaces, the same argument applies for the
existence of the invariant set as in [Martins, 2004], while the
exact number of the invariant sets may vary. More specifically, if
we cut a $2$-hole torus along a knotted trajectory, topologically
the surface will turn into a torus but with 2 boundary circles, as
illustrated in fig (\ref{cut 2-hole torus}).

\begin{figure}[!h]
\begin{center}
\includegraphics[width=4in,height=2in]{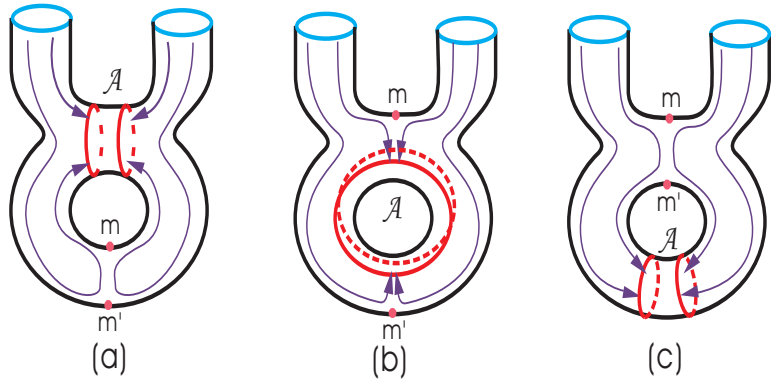}
\caption{Possible Invariant Sets in a Genus-$2$ Surface}
\label{basic attractors}
\end{center}
\end{figure}

Suppose that this knot is unstable; after the surgery, all the
dynamics are pointing outward from the two resulting boundary
circles. Since fig (\ref{cut 2-hole torus}) is essentially a
cylinder with a torus attached in the middle, from [Martins, 2004],
we know that there exists some invariant set $\mathcal{A}$. However,
the number of invariant sets differs from that of the genus-$1$
case. There can be at most three invariant sets, individually as
shown in fig (\ref{basic attractors}).

\begin{figure}[!hbp]
\begin{center}
\includegraphics[width=5in,height=2in]{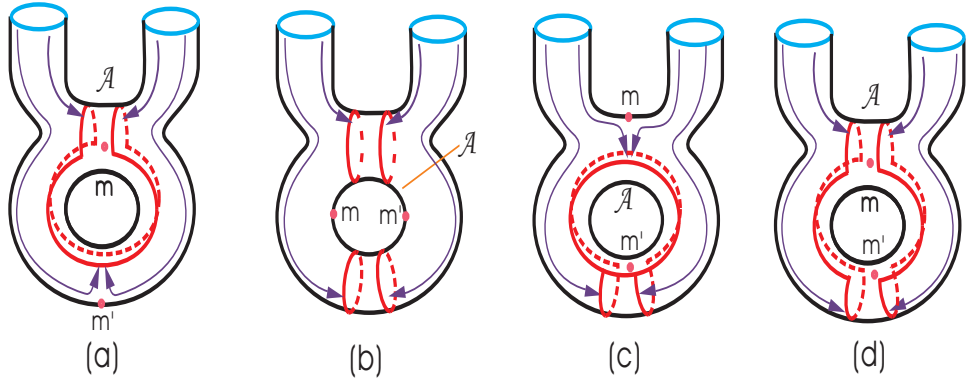}
\caption{Possible Combination of the Invariant Sets}
\label{combination of attractors}
\end{center}
\end{figure}

In this figure, $\mathcal{A}$ denotes the invariant set, while $m$
and $m'$ stand for the saddle type equilibrium points which have
$-1$ index respectively. This makes sense because
$Index(m)+Index(m')=-2$, which accounts for the correct
$Poincar\acute{e}$ characteristic for a genus-$2$ surface. Note that
the actual invariant set may be a combination of two or all of these
three individual ones (see fig (\ref{combination of attractors}) for
the illustration).

As the genus of a surface increases, the number of the invariant
sets will increase accordingly, but they are all based on the three
basic types shown in fig (\ref{basic attractors}).

\begin{lemma}
For a system situated on the genus-$p$ surface, if only one unstable
periodic orbit exists, then there can be at most $(2p-1)$ attractors
that might be knotted themselves and linked together. Moreover, a
surgery can always be performed to make them look like a combination
of the basic invariant sets shown in fig (\ref{basic attractors}).
\end{lemma}

\noindent \textbf{Proof.} We prove it by induction.

In the torus case, it is known that the attractor is a band as shown
in fig (\ref{invariant set on torus}) ([Martins, 2004]); and on a
$2$-hole torus, from the discussion above, there can be at most $3$
attractors.

Assume it is true for the genus-$p$ surface such that it has
$(2p-1)$ invariant sets at most, then by adding the genus by $1$, we
essentially introduce another hole to the manifold which will give
two more attractors at most, this proves the lemma. \qquad $\Box$

\section{Inversely Unstable Solutions and the Topology of Knotted Attractors}

Inversely unstable solutions to a dynamical system has been studied
for a long time. To be precise, we restate the main idea, which can
be found in [Levinson, 1944], for example.

Denote (\ref{periodic orbit}) by

\begin{equation} \label{inversely unstable}
\left\{
\begin{array}{l}
\dot{x}=F(x,y,t) \\
\dot{y}=G(x,y,t)
\end{array}\right.
\end{equation}

\noindent where $F$ and $G$ are both $T$-periodic in $t$.

\begin{definition}
Suppose $(a,b) \in \mathbb{Z} \times \mathbb{N}$, $b \geqslant 1$,
we shall say that a solution $z=(x,y)$ of (\ref{inversely unstable})
is $(a,b)$-periodic iff
\[
z(t+bT)=\Gamma_i^{a} \big( z(t) \big)
\]
\noindent where $\Gamma_i$ is the map between one fundamental region
and the other one next to it.
\end{definition}

Note that these solutions correspond to the trajectories that `wind
around' one of the tubes in the genus-$p$ surface $a$ times within a
$bT$ time interval before closing. If $\big( x(t),y(t) \big)$ is a
$(a,b)$-periodic solution, then the initial point $A$, $\big(
x(t_0),y(t_0)\big)$, is a fixed point of $M=P^b-\Big( \Gamma^a \big(
z(0) \big) - z(0) \Big)$. Assume $A$ is an isolated fixed point, and
let $A_0$ denote the point $\big( x(t_0)+u_0, y(t_0)+v_0 \big)$ near
$A$ in the hyperbolic upper-half plane. Applying the
$Poincar\acute{e}$ map once we have
\[
A_1=P(A_0)
\]
\noindent and $A_1$ is denoted by $\big( x(t_0)+u_1, y(t_0)+v_1
\big)$. By using a power series in $u_0$ and $v_0$ with coefficients
functions in $t$, we can express the solution trajectory of $\big(
x(t), y(t) \big)$ starting at $A_0$ by
\begin{eqnarray} \label{dynamical system}
X(t)&=& x(t)+c_1(t)u_0+c_2(t)v_0+c_3(t)u_0^2+c_4(t)u_0v_0+ \cdots \nonumber \\
Y(t)&=& y(t)+d_1(t)u_0+d_2(t)v_0+d_3(t)u_0^2+d_4(t)u_0v_0+ \cdots
\end{eqnarray}

\noindent In particular, by setting $t=t_0+T$, we have
\begin{eqnarray}
u_1 &=& au_0 + bv_0 + a_1u_0^2 + b_1u_0v_0 + \cdots \nonumber \\
v_1 &=& cu_0 + dv_0 + c_1u_0^2 + d_1u_0v_0 + \cdots
\end{eqnarray}

If we denote $\big( x(t_0)+u_0, y(t_0)+v_0 \big)$ and $\big(
x(t_0)+u_1, y(t_0)+v_1 \big)$ by $(x_0,y_0)$ and $(x_1,y_1)$
respectively, then
\begin{equation}
J(\frac{x_1,y_1}{x_0,y_0}) = J(\frac{u_1,v_1}{u_0,v_0})
\end{equation}

\noindent where $J$ is the Jacobian of the $Poincar\acute{e}$ map
for the point $(x_0,y_0)$. For very small values $u_0$ and $v_0$,
(\ref{dynamical system}) is determined by its linear terms. So the
characteristic multiplier can be determined by
\[
(a-\lambda)(d-\lambda)-bc=0.
\]

Using the notations above, we have

\begin{definition}
Given an $(a,b)$-periodic solution $\big( x(t),y(t) \big)$ of
(\ref{inversely unstable}) such that $\big( x(t_0),y(t_0) \big)$ is
an isolated fixed point of $M$, we shall say the solution $\big(
x(t),y(t) \big)$ is inversely unstable iff $\lambda_2 < -1 <
\lambda_1 <0$.
\end{definition}

In [Martins, 2004], it is shown that in the torus case, the
invariant set $\mathcal{A}$ may not be homeomorphic to a circle. We
shall now extend the ideas to higher genus surfaces. To do this we
need

\begin{definition}
A system defined on a surface S is dissipative relative to a knot K
if there is a neighbourhood, say N, of K in S such that on
$\partial(S/N)$, the vector field is pointing into N.
\end{definition}

Then we have

\begin{theorem}
Given a system defined by (\ref{periodic orbit}) on a genus-$p$
surface, which is dissipative relative to a knot K situated on this
surface as well, if there exists an inversely unstable solution
$(x_I,y_I)$ within the (knotted) attractor $\mathcal{A}_I$, then
$\mathcal{A}_I$ is not homeomorphic to the circle
$\mathbb{T}=\mathbb{R}/\mathbb{Z}$.
\end{theorem}

\begin{figure}[!h]
\begin{center}
\includegraphics[width=2.5in]{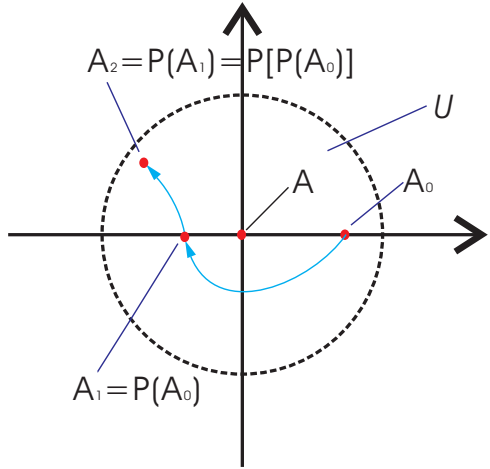}
\caption{How an Inversely Unstable Solution will Affect the
Dynamics} \label{illustration}
\end{center}
\end{figure}

\textbf{Proof.} We shall prove this theorem in a geometrical way.
Due to the dissipative nature, there exist one or more unstable
periodic orbits, and each of them is equivalent to a knot on the
surface that the system is defined on respectively. By cutting the
surface along one of these knots, we can reduce the surface genus by
$1$ while introducing two boundary circles \big(as shown in fig
(\ref{cut 2-hole torus})\big). Now gluing the two circles will
produce a tube containing an attractor $\mathcal{A}$. Assume that
there exists an inversely unstable solution in $\mathcal{A}$. Let
$A$ be a fixed point of the associated $Poincar\acute{e}$ map.
Choose a neighbourhood $U$ where $A$ is the only fixed point within
$U$. Suppose $A_0$ is a point in $U$ close to $A$ (see fig
(\ref{illustration}) for illustration). If we apply the
$Poincar\acute{e}$ map to point $A_0$, with the dynamics being
determined by the characteristic multipliers, which are $\lambda_2 <
-1 < \lambda_1 < 0$, $A_0$ will move to $A_1$, a point lies in the
other half plane with respect to $y$-axis and is much closer to the
fixed point $A$. Now apply the $Poincar\acute{e}$ map to point
$A_1$, and this time the characteristic multipliers will become $0 <
\lambda_1^2 < 1 < \lambda_2^2$ under the action of $P^2$, which
gives a directly unstable solution that moves $A_1$ to $A_2$, a
point further away in the left-half plane. With the iteration of
$Poincar\acute{e}$ map, the corresponding characteristic multipliers
will be alternatively positive and negative. However, all
neighbouring dynamics tend towards the knotted attractor by
dissipativity. In other words, within the invariant set near the
inversely unstable solution, the dynamics tend either to get close
to this trajectory or escape from it, while at the boundary, they
are pushed back by the external dissipative condition. This is why
chaotic behaviour can happen which means that $\mathcal{A}$ is not
homeomorphic to a circle. The same idea follows when there are more
than one attractor contain separate inversely unstable solutions.
\qquad $\Box$

So generally speaking, a dissipative system given by (\ref{dynamical
system}) that situated on a genus-$p$ surface can have at most $p$
topologically distinct knotted attractors; whether they are
homeomorphic to a circle individually depends on the existence of
inversely unstable solution within themselves.

It is known that any dynamical system sitting on a $2$-manifold with
$p$ genus can be represented on a sphere by cutting each handle
along a fundamental circuit which contains no equilibrium point and
filling in the dynamics within the resulting region bounded by these
curves (see [Banks, 2002]).

Conversely, we can get higher genus surface systems by performing
surgery on certain spherical ones. Specifically, given a spherical
system, irrespective of the rest of the dynamics, as long as it
contains $2$ stable equilibria, we can choose a small neighbourhood
$M_i$ $(i=1,2)$ around each of them such that they are the only
equilibrium points within each region. Glue in a dissipative region
with attractor $\mathcal{A}$ as in fig (\ref{construct torus
system}), cut this attractor open, twist it and identify the two
ends together in the appropriate way, we then obtain the desired
knot. If the attractor contains an inversely unstable solution, then
it is not homeomorphic to a circle, which means chaotic behaviour
will occur within this invariant set.

\begin{figure}[!h]
\begin{center}
\includegraphics[width=5in]{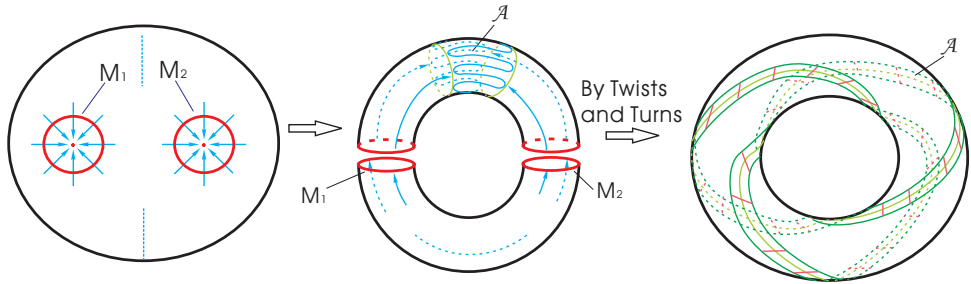}
\caption{Construct a Torus System from a Spherical One}
\label{construct torus system}
\end{center}
\end{figure}

Hence we have proved

\begin{theorem}
Any dynamical system on a genus-$p$ surface that contains a set of k
$(k\geq 1, k\in \mathbb{N})$ (knotted) dissipative attractors each
containing an inversely unstable orbit can be represented by a
system with at least 2k stable equilibrium points on a sphere.
Conversely, starting from a spherical system that contains at least
2k stable equilibria, we can construct a system on a genus-$p$
$2$-manifold that contains k knotted attractors each with chaotic
behaviour.
\end{theorem}

\noindent \textit{Remark.} An important consequence of this theorem
is that we can determine the general structure of a system with $k$
`chaotic' dissipative attractors by studying systems with $2k$
stable equilibrium points on the sphere. Of course, such a system
must have other equilibrium points so that the total index is 2, by
the $Poincar\acute{e}$ index theorem. Thus the remaining equilibrium
points must have index $2-2k$. This implies the existence of some
hyperbolic points.

\section{Examples}

In this section we show that we can obtain systems with dissipative
chaotic behaviour by choosing stable and unstable knotted orbits,
and the unstable orbit acts as the dissipative `repeller'.

\begin{figure}[!h]
\begin{center}
\includegraphics[width=3.75in,height=2in]{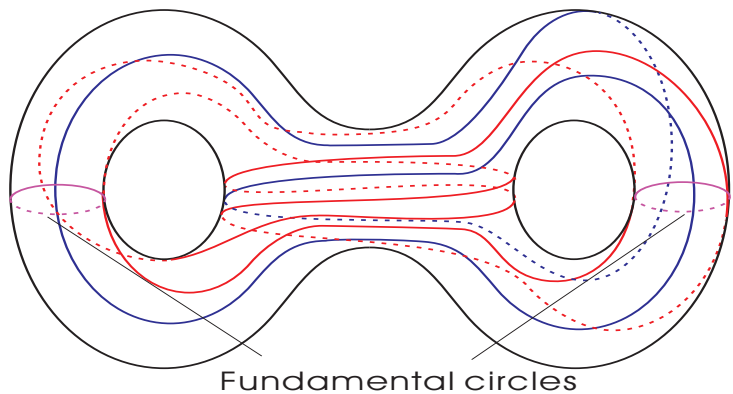}
\caption{A Surface of Genus Two Carrying Two Distinct Knot Types }
\label{2 knot types}
\end{center}
\end{figure}

In [Banks, 2002], it is shown that for a dynamical system on a
surface of genus $p$, it can carry at most $p$ distinct types of
(homotopically nontrivial) knots. For example, fig (\ref{2 knot
types}) shows the two distinct knot types that a system can have on
a genus-$2$ surface.

\begin{figure}[!h]
\begin{center}
\includegraphics[width=5in]{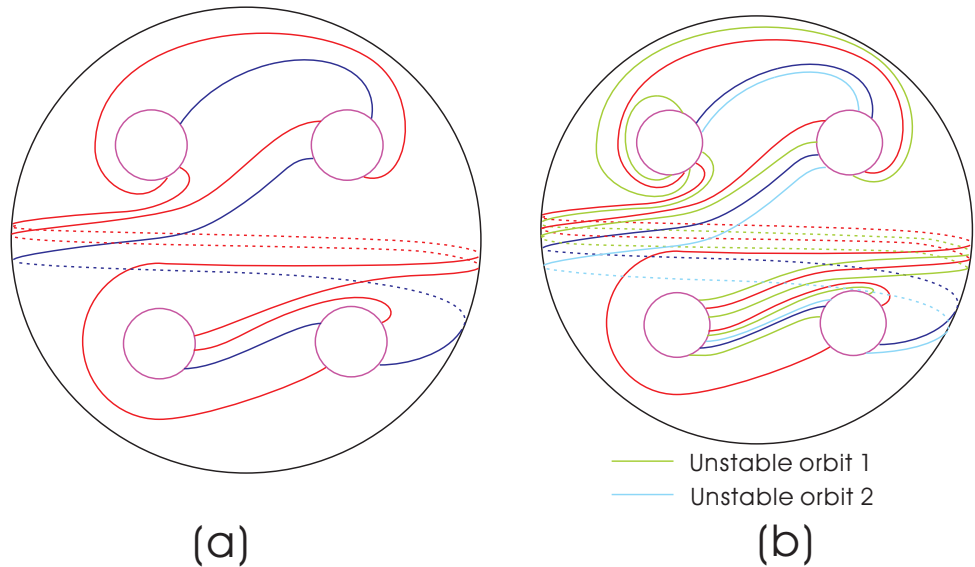}
\caption{Spherical Representation For the Attractors and The
Possible Dynamics Elsewhere } \label{spherical system}
\end{center}
\end{figure}

Assume that these two knots act as two attractors, (the existence of
chaotic behaviour will depend on whether there is an inversely
unstable solution within each attractor,) then there must exist one
or more unstable orbits due to which these two invariant sets are
generated. To find it out explicitly, we first represent the system
onto a sphere with four holes, which is achieved by cutting along
two fundamental circuits to open the handles out, as shown in fig
(\ref{spherical system}.a).

\noindent The unstable orbits therefore should bound each part of
the attractors presented on the sphere such that they can push the
dynamics toward the invariant sets and introduce possible chaotic
behaviour. Moreover, there must exist some equilibrium points to
give the correct index of a genus-$2$ surface, which is $-2$. Fig
(\ref{spherical system}.b) shows a possible solution trajectories of
two unstable orbits which satisfy the criteria discussed above.
Please note that the solution trajectories may not be unique.

\begin{figure}[!h]
\begin{center}
\includegraphics[width=4in]{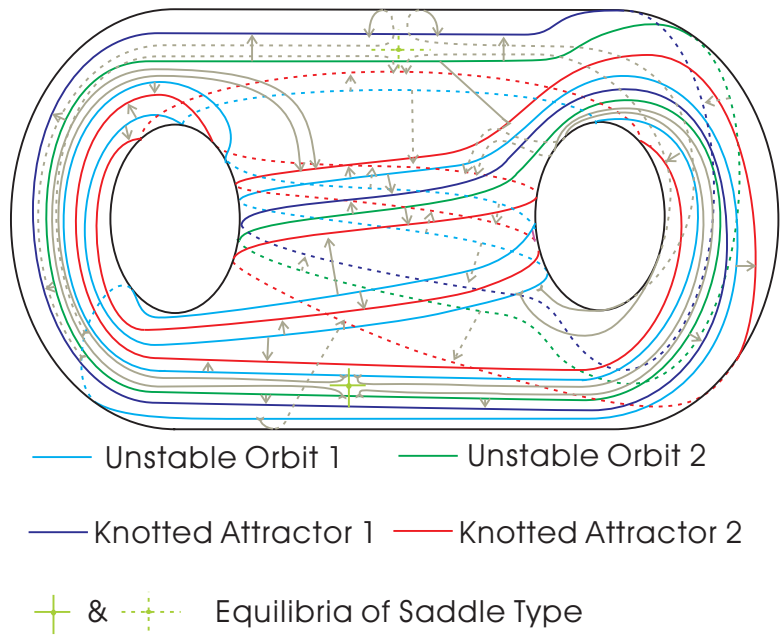}
\caption{Genus-$2$ Surface Containing 2 Knotted Attractors and the
Corresponding Dynamics}\label{example}
\end{center}
\end{figure}

Recover the original 2-manifold by gluing the corresponding boundary
circles, we eventually get a system on a genus-$2$ surface. It has
two unstable periodic cycles, which generate two knotted attractors
with distinct types, and two saddle equilibrium points which give
the correct index of $-2$ \big(See fig (\ref{example}) for an
illustration\big).

\begin{figure}[!h]
\begin{center}
\includegraphics[width=4in]{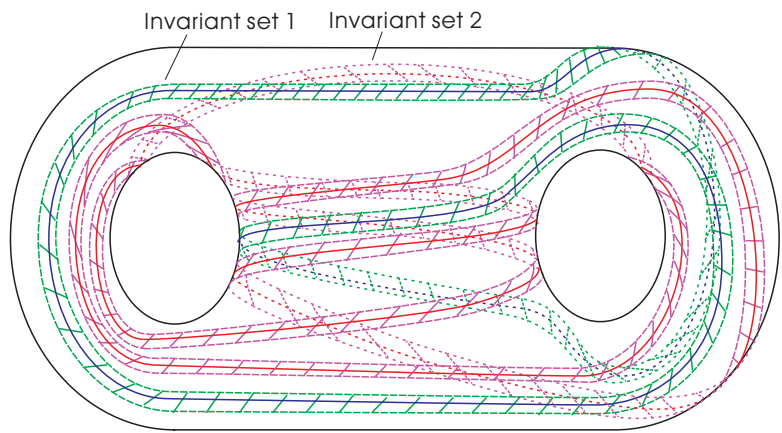}
\caption{Genus-$2$ Surface Containing 2 Invariant Sets With
Inversely Unstable Orbit In} \label{example1band}
\end{center}
\end{figure}

Moreover, as in fig (\ref{example1band}), if each invariant set
contains an inversely unstable orbit, then around each knot there
exists a band within which chaotic behaviour will occur.

\begin{figure}[!h]
\begin{center}
\includegraphics[width=4in,height=2.8in]{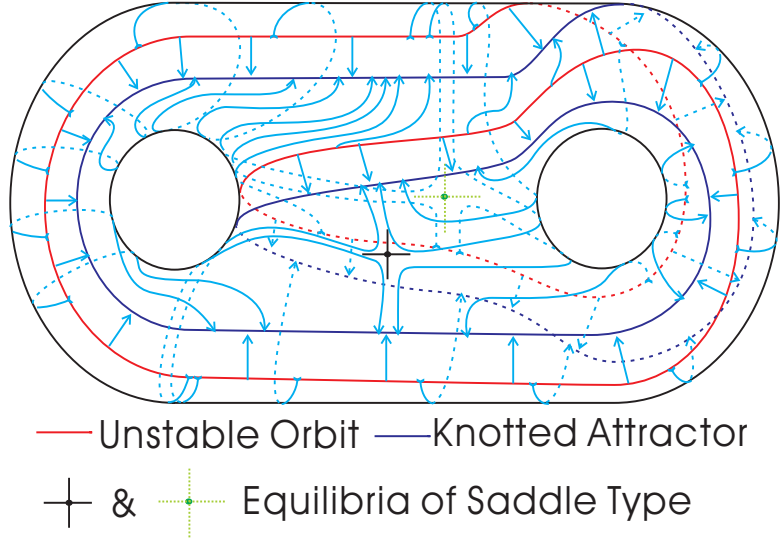}
\caption{Genus-$2$ Surface Containing 1 Knotted Attractors and the
Corresponding Dynamics} \label{example 2}
\end{center}
\end{figure}

Now if reduce the number of invariant sets by one and assume the
existence of only one unstable orbit, following the same algorithm
as above, we get one possible solution for the dynamics as in fig
(\ref{example 2}). Note that again there are two saddle equilibria
to count for the correct index.

Under the existence of inversely unstable orbit, chaotic behaviour
will occur within the invariant set \big(see fig (\ref{example 2
band})\big).

\begin{figure}[!h]
\begin{center}
\includegraphics[width=4in]{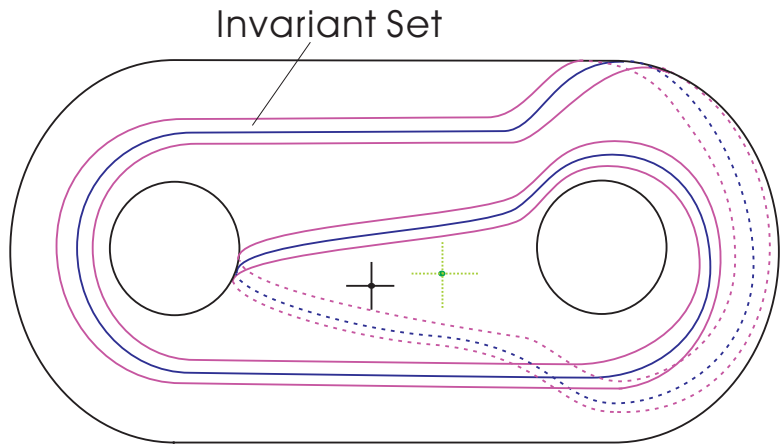}
\caption{Genus-$2$ Surface Containing 1 Knotted Attractor With an
Inversely Unstable Orbit} \label{example 2 band}
\end{center}
\end{figure}

\section{Conclusion}

We have studied dynamical systems on a genus-$p$ surface and extend
the \textit{generalized automorphic functions} (see [Banks \& Song,
2006]) to define a general form for these systems (both analytic and
non-analytic). Also we look at the topology of knotted attractors
under the existence of unstable periodic orbits and prove that for a
genus-$p$ surface with only one unstable cycle, the number of
invariant sets may vary while a maximum of $(2p-1)$ must not be
exceeded. Moreover, we extend the result in [Martins, 2004] and show
that for a higher genus $(genus > 1)$ surface, the individual
attractor is not homeomorphic to a circle if there exists an
inversely unstable solution within itself. This is purely because of
the property of inversely unstable solution which can generate a
local behaviour to make the dynamics fight against the effect of
global unstable orbit.

In the future paper, we will consider \textit{automorphic functions}
in 3-dimension which will give us systems naturally defined on
genus-$p$ solid 3-manifolds.

\end{document}